\documentclass{amsart}
\usepackage{amssymb,amscd}
\advance\hoffset by -1in
\advance\voffset by -1in

\usepackage[cp1251]{inputenc}
\usepackage[english,russian]{babel}
\newtheorem{theorem}{Theorem}
\newtheorem{lemma}{Lemma}

\newtheorem{proposition}{Proposition}
{\par\noindent{\bf Proof.}} 
{\hfill$\scriptstyle\blacksquare$}
\title{The theorems on freedom for relatively free groups with a relations}

\author{A.\,F.\,Krasnikov}
\address{Omsk State University, pr. Mira 55-A, 644077, Omsk, Russia}
\email{phomsk@mail.ru}

\begin{document}
\selectlanguage{english}
\maketitle
\markboth{A.\,F.\,Krasnikov}{The theorems on freedom}

\section*{Introduction}
The well-known result of Magnus \cite{Mg} tells us that if
$F$ --- a free group on free generators $y_1,\ldots,y_n$,
$H$ --- the subgroup of $F$, generated by $y_1,\ldots,y_{n-1}$, $1\neq r\in F$,
$R$ --- the normal subgroup of $F$ generated by $r$ then $H\cap R=1$ if and only if $r$ is not conjugate to any element of $H$.

Kolmakov \cite{Km} proved an analogous result for polynilpotent groups.

In \S1 is contained the definition of Fox's derivations and various considerations in the group algebra of the free group. The prove the theorem on freedom for relatively free groups with a single relation occupies \S2:
\begin{theorem}\label{tm_1}
Suppose $F$ is a free group on free generators $y_1,\ldots,y_n\,(n>2)$,
$H$ --- the subgroup of $F$, generated by $y_1,\ldots,y_{n-1}$, $N_{11}$  ---  a normal subgroup of $F$, $F/N_{11}$  --- an orderable and relatively free group,
\begin{eqnarray}\label{tm_1_f_1}
N_{11} \geqslant \ldots \geqslant N_{1,m_1+1}=N_{21} \geqslant \ldots \geqslant N_{s,m_s+1},
\end{eqnarray}
where $N_{kl}$ --- the $l$-th term of the lower central series of $N_{k1}$.
Let $r$ be an element of $N_{1i}\backslash N_{1,i+1}\,(i\leqslant m_1)$, $R$ --- the normal subgroup of $F$ generated by $r$.
If (and only if) $r$ is not conjugate to any element of $HN_{1,i+1}$ then $H\cap RN_{kl} = H\cap N_{kl}$ for each term $N_{kl}$ of series {\rm (\ref{tm_1_f_1})}.
\end{theorem}

If $G$ is a free polynilpotent group then $G$ --- an orderable group \cite{Sh}. Hence from the theorem \ref{tm_1_f_1} we have the theorem on freedom for polynilpotent groups with a single relation \cite{Km}.

Romanovskii \cite{Rm2} proved that if $F$ --- a free group on free generators $y_1,\ldots,y_n$,
$r_1,\ldots,r_m$ --- an elements of $F$ $(m<n)$, $R$ --- the normal subgroup of $F$ generated by $r_1,\ldots,r_m$,
then there exists a subset $y_{j_1},\ldots,y_{j_p}$ $(p\geqslant n-m)$ of $y_1,\ldots,y_n$
such that $H\cap R=1$, where $H$ is the subgroup of $F$ generated by $y_{j_1},\ldots,y_{j_p}$.

Romanovskii proved also an analogous results for solvable groups \cite{Rm2}, pro-$p$-groups \cite{Rm3}
and for products of groups \cite{Rm4}.

The prove the theorem on freedom for relatively free groups with a relations (generalized Freiheitssatz) occupies \S3:
\begin{theorem}\label{tm_2}
Suppose $F$ is a free group on free generators $y_1,\ldots,y_n$,
$N_{11}$  ---  a normal subgroup of $F$, $F/N_{11}$  --- a soluble, right-ordered and relatively free group,
\begin{eqnarray}\label{tm_2_f_1}
N_{11} \geqslant \ldots \geqslant N_{1,m_1+1}=N_{21} \geqslant \ldots \geqslant N_{s,m_s+1},
\end{eqnarray}
where $N_{kl}$ --- the $l$-th term of the lower central series of $N_{k1}$.
Let $r_1,\ldots,r_m$ be an elements of $N_{11}$ $(m<n)$, $R$ --- the normal subgroup of $F$ generated by $r_1,\ldots,r_m$. Then there exists a subset $y_{j_1},\ldots,y_{j_p}$ $(p\geqslant n-m)$ of $y_1,\ldots,y_n$
such that $H\cap RN_{kl} = H\cap N_{kl}$ for each term $N_{kl}$ of series {\rm (\ref{tm_2_f_1})}, where $H$ is the subgroup of $F$ generated by $y_{j_1},\ldots,y_{j_p}$.
\end{theorem}

From the theorem \ref{tm_2} we have the theorem on freedom for polynilpotent groups with a relations.

\section{Certain properties of the Fox derivations}
With any group $G$ there is associated its group ring ${\bf Z}[G]$. An element of ${\bf Z}[G]$ is a sum
$\sum a_g g,\,g$ ranging over the elements of $G$, where the integer $a_g$ is equal to zero for all but a finite number of $g$.
Addition and multiplication in ${\bf Z}[G]$ are defined by $\sum a_g g+\sum b_g g=\sum (a_g+b_g) g$
and $(\sum a_g g)\cdot (\sum b_g g)=\sum(\sum_h a_{gh^{-1}}b_h) g$. The element $a$ of ${\bf Z}$ is identified with the element $a\cdot 1$ of ${\bf Z}[G]$ and the element $g$ of $G$ is identified with the element $1\cdot g$ of ${\bf Z}[G]$, so that ${\bf Z}$ and $G$ are to be regarded as subsets of ${\bf Z}[G]$.
By a derivation in a group ring ${\bf Z}[G]$ will be meant any mapping $\partial$ of ${\bf Z}[G]$ into itself which satisfies $\partial(u+v)=\partial(u)+\partial(v)$, $\partial(uv)=\partial(u)v+\varepsilon (u)\partial(v)$ for any $u,\,v\in {\bf Z}[G]$, where $\varepsilon:\,{\bf Z}[G]\to {\bf Z}$ is the natural augmentation.

Let $F$ be a free group on free generators $\{g_j | j\in J\}$, $N$ --- a normal subgroup of $F$. We denote by $D_j$ $(j\in J)$ the Fox derivations of the group ring ${\bf Z}[F]$. They are uniquely defined by the conditions $D_j(g_j)=1,\,D_k(g_j)=0\,\mbox{for}\,k\neq j$. It is easy to deduce the following relations:
\begin{gather}
D_j(f^{-1}nf)\equiv D_j(n)f\mod{{\bf Z}[F]\cdot (N-1)};\label{fm00}\\
D_j(f^{-1})=-D_j(f)f^{-1};\notag\\
u-\varepsilon(u)= \sum_{j\in J} (g_j-1)D_j(u);\notag
\end{gather}
where $u\in {\bf Z}(F),\,f\in F,\,n\in N$.

Let $G$ be a group $x,\,y\in G$; $X,\,Y$ are subgroups in $G$. Then $[x,y]=x^{-1}y^{-1}xy$, $x^y=y^{-1}xy$;
$[X,Y]$ denotes the subgroup generated by all commutators of the type $[x,y]$ and $XY$ denotes the subgroup generated by
$\{xy\mid x\in X,\,y\in Y\}$.
If $A$ is subgroup in $G$, then $A^G$ denotes the normal subgroup generated by $A$ in $G$.
By $\gamma_k (G)$ we denote the $k$-th term of the lower central series of a group $G$, i.e.
$\gamma_1(G)=G$, $\gamma_{i+1}(G)=[\gamma_i(G),G]$.

\begin{lemma}\label{lm_3}
Let $X$ be a free group on free generators $\{x_j \mid j\in J\}$, $X_n$ --- $n$-th term of the lower central series of $X$,
$\mathfrak{X}$ --- the fundamental ideal of ${\bf Z}[X]$, $v\in X$. Then $v\in X_n\setminus X_{n+1}$ if and only if $D_j(v)\in \mathfrak{X}^{n-1} \, (j\in J)$
and there exists $j_0\in J$ such that $D_{j_0}(v)\in \mathfrak{X}^{n-1}\setminus \mathfrak{X}^n$.
\end{lemma}
\begin{proof}
It is known \cite{Fx} that quotient ring $\mathfrak{X}^k/\mathfrak{X}^{k+1}$ has an explicit basis $(x_{j_1}-1)\ldots (x_{j_k}-1)+\mathfrak{X}^{k+1}$ and $v\in X_n$  if and only if $v-1\in \mathfrak{X}^n$.
Now we reach the conclusion by applying
$v-1=\sum_{j\in J} (x_j-1)D_j(v)$.
\end{proof}

Using the induction on the length of the word $v$ it is easy to proof
\begin{lemma}\label{lm_5}
Let $G$ be a free group on free generators $\{g_j \mid j\in J\}$, $\{D_j | j\in J\}$ --- the Fox derivations of the group ring ${\bf Z}[G]$, $H$ --- a subgroup of $G$ with a free set $\{x_i \mid i\in I\}$ of generators,
$\{\partial_i \mid i\in I\}$ --- the Fox derivations of the group ring ${\bf Z}[H]$, $v\in H$.
Then
\begin{eqnarray}\label{lm_5_f_1}
D_j(v)=\sum_k D_j(x_k)\partial_k(v).
\end{eqnarray}
\end{lemma}

\begin{lemma}\label{lm_6}
Let $X$ be a free group on free generators $\{x_j \mid j\in {\bf N}\}$,
$\mathfrak{X}$ --- the fundamental ideal of ${\bf Z}[X]$, $v\in X$. If
$D_1(v)\notin \mathfrak{X}^{j-1}$ then $D_1([v,x_2])\notin \mathfrak{X}^j$.
\end{lemma}
\begin{proof}
Since $D_1([v,x_2])=D_1(v)(x_2 - 1) + D_1(v)(1-v^{-1}{x_2}^{-1}vx_2)$ and
$1-v^{-1}{x_2}^{-1}vx_2\in \mathfrak{X}^2$ \cite{Fx} it follows that $D_1([v,x_2])\notin \mathfrak{X}^j$.
\end{proof}

\begin{lemma}\label{lm_7}
Let $G$ be a free group on free generators $\{g_j \mid j\in J\}$, $\{D_j | j\in J\}$ --- the Fox derivations of the group ring ${\bf Z}[G]$, $N$ --- a normal subgroup of $G$, $v\in N$, $S$ --- a Schreier system of representatives of $G$ by $N$, $u\rightarrow \bar{u}$ --- a Schreier coset representative function for $G\mod{N}$, $\{x_i \mid i\in I\}$--- a free set of generators of $N$ such that $\{x_i \mid i\in I\}\subseteq \{sg_j\overline{sg_j}^{\,-1} \mid s\in S,\, j\in J\}$.
If $x_i=sg_j\overline{sg_j}^{\,-1}$ then $D_j(v)=\overline{sg_j}^{\,-1} \partial_i(v)+V$,
where $V$ --- a linear combination of an elements of the form $g t$, $t\in {\bf Z}[N]$, $g\in S$, $g\not\equiv \overline{sg_j}^{\,-1}\mod{N}$.
\end{lemma}
\begin{proof}
We have $D_j(x_i)= D_j(s){s}^{-1}x_i+
\overline{sg_j}^{\,-1}-D_j(\overline{sg_j})\overline{sg_j}^{\,-1}$.\\
If $x_k=tg_l\overline{tg_l}^{\,-1}$ $(t\in S,\,k\neq i)$ then
$D_j(x_k)= D_j(t){t}^{-1}x_k + \alpha\cdot\overline{t g_l}^{\,-1}-
D_j(\overline{t g_l}){\overline{t g_l}}^{-1}$ where $\alpha=1$ for $l=j$
and $\alpha=0$ for $l\neq j$.
If $l=j$ then $\overline{tg_l}\neq \overline{sg_j}$.

By Lemma \ref{lm_5} $D_j(v)=\sum_k D_j(x_k)\partial_k(v)$, hence we need to prove that from $h\in S$ follows that $D_j(h){h}^{-1}$ --- a linear combination of an elements of the form $g t$, $t\in {\bf Z}[N]$, $g\in S$, $g\not\equiv \overline{sg_j}^{\,-1}$.
Let $h=ug^\varepsilon_j u_1$ where $\varepsilon = \pm 1$ and $D_j(u_1)=0$.
Then $D_j(h){h}^{-1}=D_j(u){u}^{-1}+D_j(g^\varepsilon_j)(ug^\varepsilon_j)^{-1},$
i.e. $D_j(h){h}^{-1}$ --- a linear combination of an elements of the form $D_j({g_j}^\varepsilon)(ug^\varepsilon_j)^{-1}$, $ug^\varepsilon_j\in S$.

We will show that $\overline{sg_j}^{\,-1}\not\equiv \pm D_j({g_j}^\varepsilon)({u{g_j}^\varepsilon})^{-1}\mod{N}.$
Assume the contrary. We have $\overline{sg_j}^{\,-1}\equiv (u{g_j})^{-1}\mod{N}$ for $\varepsilon = 1$.
Hence $s=u$ and $\overline{sg_j}=sg_j$, a contradiction.
We have $\overline{sg_j}^{\,-1}\equiv {u}^{-1}\mod{N}$ for $\varepsilon = -1$.
Hence $s=ug_j^{-1}$ and $sg_j\overline{sg_j}^{\,-1}=1$, a contradiction.
The lemma is proved.
\end{proof}

\begin{lemma}\label{tm_3}\cite{Kr}
Let $F$ be a free group on free generators $\{g_j \mid j\in J\}$, $K\subseteq J$, $F_K$
 --- the subgroup of $F$ generated by $\{g_j\mid j\in K\}$; $v\in F$,
$N$--- a normal subgroup of $F$; $D_k (k\in J)$ --- the Fox derivations of the group ring ${\bf Z}[F]$. Then $D_k(v)\equiv 0\mod{{\bf Z}[F]\cdot(N-1)}, k\in J\setminus K$ if and only if $v\in F_K(F_K\cap N)^F\mod{[N,N]}$.
\end{lemma}

\begin{lemma}\label{gr1_lm2}
Let $X$ be a free group on free generators $\{x_j \mid j\in J\}$,
$\mathfrak{X}$ --- the fundamental ideal of ${\bf Z}[X]$, $v\in X$, $K\subseteq J$, $X_K$ --- a subgroup of $X$ generated by $\{x_j\mid j\in K\}$. Then from
\begin{eqnarray*}
D_k(v)\equiv 0\mod{\mathfrak{X}^n},\,k\in J\setminus K;\,D_k(v)\in{\bf Z}[X_K]\mod{\mathfrak{X}^n},\,k\in  K
\end{eqnarray*}
follows that $v$ modulo $\gamma_{n+1}(X)$ is an element of $X_K$.
\end{lemma}
\begin{proof}
We consider a homomorphism $\varphi\text{: }{\bf Z}[X]\rightarrow {\bf Z}[X]$ which is defined by the mapping
$x_j\rightarrow x_j$ for $j\in K$, $x_j\rightarrow 1$ for $j\in J\setminus K$ and denote $\varphi(v)$ by $\bar{v}$.
It is clear that $\bar{v}\in F_K$ and $D_k(v)\equiv D_k(\bar{v})\mod{\mathfrak{X}^n}\,(k\in J)$.
Since $D_k(v\bar{v}^{-1})=D_k(v)\bar{v}^{-1}-D_k(\bar{v})\bar{v}^{-1}$ it follows that
$D_k(v\bar{v}^{-1})\equiv 0\mod{\mathfrak{X}^n}\,(k\in J)$ and we have $v\bar{v}^{-1}-1\in \mathfrak{X}^{n+1}$.
Hence $v\bar{v}^{-1}\in \gamma_{n+1}(X)$ .
\end{proof}

\section{A theorem on freedom for relatively free groups with a single relation}

Let $F$ be a free group on free generators $\{g_j \mid j\in J\}$, $K\subseteq J$,  $H$ --- a subgroup of $F$ generated by $\{g_j\mid j\in K\}$; $N=N_1 \geqslant \ldots \geqslant N_t \geqslant \ldots$ --- a descending series of normal subgroups of $F$ with abelian torsion free factors such that $[N_i,N_j\,]\leqslant N_{i+j}$, $H_i=H\cap N_i$.
We denote the ideal in ${\bf Z}[N]$ generated by all products of the form $(N_{i_1} - 1)\cdots (N_{i_t} - 1)$ $(i_1 + \cdots + i_t\geqslant i)$ by $\Delta_i$; the ideal in ${\bf Z}[H_1]$ generated by all products of the form $(H_{i_1} - 1)\cdots (H_{i_t} - 1)$ $(i_1 + \cdots + i_t\geqslant i)$ by $\Delta_i^\prime$.
Define also that $\Delta_0={\bf Z}[N]$, $\Delta_0^\prime={\bf Z}[H_1]$.

A coset of $F\mod{N}$ we call $\alpha$-coset if there exist an elements of $H$ in it and $\beta$-coset if an elements of $H$ are not in it.
We call the {\it length} of a $\alpha$-coset the length of the shortest word of $H$ in it;
the {\it length} of a $\beta$-coset the length of the shortest word in it and
shall define Schreier representatives inductively, using the {\it length} of a coset.

Choose the empty word as a representative of $N$. If $L$ is a $\alpha$-coset of {\it length} one, choose any word of length one
in $L\cap H$ as it representative. Assuming we have chosen representatives for all $\alpha$-cosets of {\it length} less then $r$, i.e. we have a Schreier coset function $u\rightarrow \bar{u}$ for all $\alpha$-cosets of {\it length} less then $r$.
If $L$ is a $\alpha$-coset of {\it length} $r$ and $z_1\ldots z_r$ $(z_m\in \{g_j^{\pm 1} \mid j\in K\})$ is a word in $L$,
we choose $\overline{z_1\ldots z_{r-1}}z_r$ as a representative of $L$. Denote by $S_\alpha$ the system of representatives of $\alpha$-cosets.

If $L$ is a $\beta$-coset of {\it length} one, choose any word of length one in $L$ as it representative. Assuming we have chosen representatives for all $\alpha$-cosets and all $\beta$-cosets of {\it length} less then $r$. If $L$ is a $\beta$-coset of {\it length} $r$ and $z_1\ldots z_r$ $(z_m\in \{g_j^{\pm 1} \mid j\in J\})$ is a word in $L$, we choose $\overline{z_1\ldots z_{r-1}}z_r$ as a representative of $L$. Denote by $S_\beta$ the system of representatives of $\beta$-cosets and by $S$ --- a system of representatives $S\alpha\cup S_\beta$.

Let $G$ be a finitely generated subgroup of $N$, $G_i=G\cap N_i$, $H_i^\prime=H\cap G_i$.
We denote the ideal in ${\bf Z}[G]$ generated by all products of the form $(G_{i_1} - 1)\cdots (G_{i_t} - 1)$ $(i_1 + \cdots + i_t\geqslant i)$ by $\delta_i$; the ideal in ${\bf Z}[H_1^\prime]$ generated by all products of the form $(H_{i_1}^\prime - 1)\cdots (H_{i_t}^\prime - 1)$ $(i_1 + \cdots + i_t\geqslant i)$ by $\delta_i^\prime$.
Define also that $\delta_0={\bf Z}[G]$, $\delta_0^\prime={\bf Z}[H_1^\prime]$.

We choose coordinated bases $M_k=a_{k,1},\ldots,a_{k,I_k}$, $b_{k,1},\ldots,b_{k,J_k}$ for $G_k$ modulo $G_{k+1}$
such that the elements of the coordinated bases for $H_k^\prime$ modulo $H_{k+1}^\prime$
will be powers modulo $G_{k+1}$ of elements $a_{k,1},\ldots,a_{k,I_k}$.
We denote $i$-th term of series $M_k$ by $c_{ki}$ and
define by $C_{kt}$ the subgroup of $G$, generated by $\{c_{ki} \mid i\geqslant t\}$ and $G_{k+1}$.

We introduce on ${\bf N}\times {\bf N}$ the lexicographic ordering $(s,m)<(k,j)$ if $s<k$ or $s=k$ and $m<j$ and
name the {\it monomial} a product
\begin{eqnarray}\label{alg2_0}
(c_{j_1 k_1}-1)\ldots (c_{j_\nu k_\nu}-1),
\end{eqnarray}
where $(j_1,k_1)\leqslant\ldots\leqslant (j_\nu,k_\nu),\,\nu\geqslant 0$.
Define the {\it weight} of the monomial (\ref{alg2_0}) as $j_1+\ldots + j_\nu$.
Define also that $1$ is a monomial of a weight $0$.
Order monomials by their weights and
order monomials of equal weight by taking monomials of greater
length to be greater than monomials of smaller length, and using
lexicographic order for monomials of equal length (from left to right).

It can be verified directly that
\begin{gather}
n(a-1)\equiv (a^n-1)\mod {\Delta_{k+1}},\text{ }a \in N_k\setminus N_{k+1};\label{fm01}\\
(a-1)(b-1)=(b-1)(a-1)+ba([a,b]-1),\text{ }a,\,b \in N.\label{fm02}
\end{gather}

Now we state some known results and include short proofs.

1) Monomials of a weight $i$, which are not in ${\bf Z}[G](C_{lt}-1)$, form a linear basis over ${\bf Z}$ of $\delta_i$ modulo $\delta_{i+1} + {\bf Z}[G](C_{lt}-1),\,(l,t)\neq (1,1)$.
\begin{proof} The fact that monomials of a weight $i$ over ${\bf Z}$ generate the ideal $\delta_i$ modulo $\delta_{i+1}$ can be proved by a collection process based on the formula (\ref{fm02}).
Note that one can prove that every element of $\delta_i$ has unique representation modulo $\delta_{i+1} + {\bf Z}[G](C_{lt}-1)$ as a linear combination with coefficient from ${\bf Z}$ of monomials of a weight $i$ which are not in ${\bf Z}[G](C_{lt}-1),\,(l,t)\neq (1,1)$. This is true if $i=0$.
Assume inductively that this statement is true for $i\leqslant k-1$. We need to show that this statement is true for $i=k$
and any $(l,t)\neq (1,1)$. Assume the contrary and choose a minimal $(l_0,t_0)\neq (1,1)$ which contradicts the statement.

Let $v$ be a linear combination with non-zero coefficient from ${\bf Z}$ of monomials of a weight $k$ such that $v\in \delta_{k+1} + {\bf Z}[G](C_{l_0t_0}-1)$. Then $v=u(c_{lt}-1)$, where $c_{lt}=\max\{c_{ij}\mid c_{ij}<c_{l_0t_0}\}$ and
$u$ is a linear combination with non-zero coefficient from ${\bf Z}$ of monomials of a weight $k-l$ which are not in ${\bf Z}[G](C_{l_0t_0}-1)$.
Since ${\bf Z}[G/C_{l_0t_0}]$ has no zero-divisors, we obtain that $u$ --- a linear combination
with coefficient from ${\bf Z}$ of monomials of a weight $k-l$ and $u\in \delta_{k-l+1} + {\bf Z}[G](C_{l_0t_0}-1)$. A contradiction.
The statement is proved.
\end{proof}

2) Let $M_1$, $M_2$ be monomials of the form (\ref{alg2_0}). We denote by $M_1\circ M_2$ a monomial which is produced from the product $M_1M_2$ by corresponding formal substitution of factors $c_{ij}-1$.
Using the collection process, one can easily prove that if $M_1,\,M_2$ are monomials of a weight $i,\,j$ respectively then $M_1\circ M_2$ is a maximal monomial involved in the representation of $M_1M_2$ modulo $\delta_{i+j+1}$ as a linear combination with coefficient from ${\bf Z}$ of monomials of the weight $i+j$.

Let $u_1,\,u_2$ be monomials of the form (\ref{alg2_0}) of a weight $i$, $v_1,\,v_2$ --- monomials of the form (\ref{alg2_0}) of a weight $j$. It can be verified directly that if $u_1<u_2$, $v_1\leqslant v_2$ then $u_1\circ v_1<u_2\circ v_2$.
Hence, if $u\in \Delta_i\setminus \Delta_{i+1}+{\bf Z}[F](N_m-1)$, $v\in \Delta_j\setminus \Delta_{j+1}+{\bf Z}[F](N_m-1)$ then
$uv\in \Delta_{i+j}\setminus \Delta_{i+j+1}+{\bf Z}[F](N_m-1)$.
Since $F/N$ is a right-ordered group, it is not hard to verify that
if $u\in S\Delta_i\setminus S\Delta_{i+1}+{\bf Z}[F](N_m-1)$, $v\in S\Delta_j\setminus S\Delta_{j+1}+{\bf Z}[F](N_m-1)$ then
$uv\in S\Delta_{i+j}\setminus S\Delta_{i+j+1}+{\bf Z}[F](N_m-1)$.

Let $\phi$ be a natural homomorphism ${\bf Z}[F]\to {\bf Z}[F/N_m]$, $u,\,v\in {\bf Z}[F]$. Define a function $\psi$ on ${\bf Z}[F/N_m]$ by $\psi(\phi(u))=j$ if $u\in S\Delta_j\setminus S\Delta_{j+1}\mod{{\bf Z}[F](N_m-1)}$, $\psi(0)=\infty$.
Since $\psi(\phi(uv))=\psi(\phi(u))+\psi(\phi(v))$ and $\psi(\phi(u+v))\geqslant \min\{\psi(\phi(u)),\psi(\phi(v))\}$
it follows that $\psi$ is a valuation.

3) Let $u\in {\bf Z}[H_1]$. Then $u\in {\bf Z}[H_1]\cap \delta_i$ if and only if $u\in \delta_i^\prime$.
\begin{proof}
The inclusion $\delta_i^\prime\leqslant {\bf Z}[H_1]\cap \delta_i$ is obvious.
Let $u\in {\bf Z}[H_1]\cap \delta_i$.
The elements of the coordinated bases for $H_k^\prime$ modulo $H_{k+1}^\prime$
will be powers modulo $G_{k+1}$ of elements $a_{k,1},\ldots,a_{k,I_k}$.
By (\ref{fm01}) we obtain that monomials of a weight $l$ of an element $u$ in the first representation will be
coordinated with the monomials of a weight $l$ of $u$ in the second representation
The required inclusion $u\in \delta_i^\prime$ follows from this.
Thus we have
\begin{eqnarray}\label{fm04}
{\bf Z}[H]\cap \Delta_i=\Delta_i^\prime.
\end{eqnarray}
\end{proof}

\begin{lemma}\label{lm_1}
Let $F$ be a free group on free generators $\{g_j \mid j\in J\}$, $K\subseteq J$,  $H$ --- a subgroup of $F$ generated by $\{g_j\mid j\in K\}$; $N=N_1 \geqslant \ldots \geqslant N_t \geqslant \ldots$ --- a descending series of normal subgroups of $F$ with abelian torsion free factors such that $[N_i,N_j\,]\leqslant N_{i+j}$, $S=S_\alpha\cup S_\beta$ --- a system of representatives of $F$ by $N$, $F/N$ --- right-ordered group. Let also that
$r \in S_\alpha\Delta_{t-1}\mod{{\bf Z}[F]\cdot\Delta_t}$, $r \not\in {\bf Z}[F]\cdot\Delta_t$, $w \in S_\alpha\Delta_{l-t}\mod{{\bf Z}[F]\cdot\Delta_{l-t+1}}$, $w \not\in {\bf Z}[F]\cdot\Delta_{l-t+1}$. Then from $rw\in {\bf Z}[H]\mod{{\bf Z}[F]\cdot\Delta_l}$
follows that there exists $M\in {\bf N}$ such that
$Mr\in {\bf Z}[H]\mod{{\bf Z}[F]\cdot\Delta_t}$.
\end{lemma}
\begin{proof}
Let $G$ be a finitely generated subgroup of $N$, $M_k=a_{k,1},\ldots,a_{k,I_k}$, $b_{k,1},\ldots,b_{k,J_k}$ --- coordinated bases for $G_k$ modulo $G_{k+1}$ where $a_{k,1},\ldots,a_{k,I_k}$ such that the elements of the coordinated bases for $H_k^\prime$ modulo $H_{k+1}^\prime$ will be powers modulo $G_{k+1}$ of elements $a_{k,1},\ldots,a_{k,I_k}$.

We define function $d_G$ on ${\bf Z}[G]$ by
$d_G(c_{k i}-1)=1 \text{ for }i>I_k\text{ and }d_b(c_{k i}-1)=0\text{ for }i\leqslant I_k $;
$d_G(v)=\sum_{i=1}^\nu d_b(c_{j_i k_i}-1)$ if $v$ is a monomial (\ref{alg2_0});
$d_G(u)=\max \{d_G(u_i)\mid i=1,\ldots,k\}$ if $u\in \delta_j\setminus \delta_{j+1}$ has a unique representation modulo $\delta_{j+1}$ as a linear combination with coefficient from ${\bf Z}$ of monomials $u_1,\ldots,u_k$ of a weight $j$.

It is not hard to verify that if $G_1,\,G_2$ are a finitely generated subgroups of $N$, $G_1\leqslant G_2$, $u\in {\bf Z}[G_1]$
then $d_{G_1}(u)\geqslant d_{G_2}(u)$. Thus we may, without ambiguity, define function $d_b$ on ${\bf Z}[N]$ by $d_b(u)=d_G(u)$ where
$G$ --- finitely generated subgroup of $N$ such that $u\in {\bf Z}[G]$ and $d_G(u)= d_{G^\prime}(u)$ for any $G^\prime\geqslant G$.
Using (\ref{fm01}), (\ref{fm02}) one can easily proof that $d_b(uv)=d_b(u)+d_b(v)$ where $u,\,v\in {\bf Z}[N]$ and
can check directly that if $f\in H$, $u\in {\bf Z}[N]$ then $d_b(u)=d_b(u^f)$.

Function $d_b$ may be extended to a function defined on the ${\bf Z}[F]$ by the equation $d_b(u)=\max\{d_b(u_1),\ldots,d_b(u_x)\}$ where $u=g_1 u_1+\ldots +g_x u_x$; $g_1,\ldots,g_x\in S$ $(g_i\neq g_j\mbox{ for }i\neq j)$; $u_1,\ldots,u_x\in {\bf Z}[N]$.

Let $r=f_1 A_1+\ldots +f_x A_x$; $f_1,\ldots,f_x\in S_\alpha$ $(f_i\neq f_j\mbox{ for }i\neq j)$; $A_1,\ldots, A_x\in {\bf Z}[N]$;
$w=g_1 B_1+\ldots +g_y B_y$; $g_1,\ldots,g_y\in S_\alpha$ $(g_i\neq g_j\mbox{ for }i\neq j)$; $B_1,\ldots, B_y\in {\bf Z}[N]$;
$A=\{f_p\mid d_b(A_p)=d_b(r)\}$, $B=\{g_k\mid d_b(B_k)=d_b(w)\}$.
Since $F/N$ is a right-ordered group, it follows that there exist
$f_{p_0}\in A$, $g_{k_0}\in B$ such that $f_{p_0}g_{k_0}\not\equiv f_pg_k\mod{N}$
for $(p_0,k_0)\neq (p,k)$, $f_p\in A$, $g_k\in B$. We have $0=d_b(rw)=d_b(A_{p_0}^{g_{k_0}}B_{k_0})=d_b(r)+d_b(w)$.
Thus $d_b(r)=0$ whence follows the  existence
 $M\in {\bf N}$ such that
$Mr\in {\bf Z}[H]\mod{{\bf Z}[F]\cdot\Delta_t}$. The statement is proved.
\end{proof}

Let $F$ be a free group, $N$, $R$ --- a normal subgroups of $F$, $N\geqslant R$, $F/N$ --- right-ordered group,
$N=N_{11} \geqslant \ldots \geqslant N_{1,m_1+1}=N_{21} \geqslant \ldots \geqslant N_{s,m_s+1}$,
where $N_{kl}$ --- the $l$-th term of the lower central series of $N_{k1}$.

We define subgroups $\sqrt{RN_{kl}}$ of $F$ inductively as follows: $\sqrt{RN_{11}}=N$;
$\sqrt{RN_{kl}}=\{x\in \sqrt{RN_{k1}}\mid x^n\in \gamma_l (\sqrt{RN_{k1}})R\text{ for some }n\neq 0\}$; $\sqrt{RN_{k+1,1}}\,=\sqrt{RN_{k,m_k+1}}$.

Since $\sqrt{RN_{k1}}/(\gamma_l (\sqrt{RN_{k1}})R)$ is a nilpotent group it follows that
$\sqrt{RN_{kl}}$ --- a subgroup of $\sqrt{RN_{k1}}$. Thus $\sqrt{RN_{kl}}$ --- a normal subgroup of $F$.

Since $F/N_{11}$ is a right-ordered group and group $N_{11}/\sqrt{RN_{kl}}$ has a normal series with abelian torsion free factors
it follows that $F/\sqrt{RN_{kl}}$ is a right-ordered group.

Since $\sqrt{RN_{k1}}/\sqrt{RN_{k,p+m}}$ is a nilpotent torsion free group
it follows that if $x,\, y$ --- an elements of $\sqrt{RN_{k1}}/\sqrt{RN_{k,p+m}}$ and $x^t y^n=y^n x^t\,(t,\,n\neq 0)$ then $x y=y x$. Thus $[\sqrt{RN_{kp}}\,,\sqrt{RN_{km}}\,]\leqslant \sqrt{RN_{k,p+m}}$.

\begin{proposition}\label{prp_1}
Suppose $F$ is a free group on free generators $y_1,\ldots,y_n$, $n\geqslant 3$, $1\neq N_{11}$ --- a normal subgroup of $F$,
$F/N_{11}$ --- right-ordered and relatively free group,
\begin{eqnarray}\label{tm4_0}
N_{11} \geqslant \ldots \geqslant N_{1,m_1+1}=N_{21} \geqslant \ldots \geqslant N_{s,m_s+1},
\end{eqnarray}
where $N_{kl}$ --- the $l$-th term of the lower central series of $N_{k1}$.
Let $R$ be a normal subgroup of $F$, $R\leqslant N$; $H$ --- the subgroup of $F$, generated by $y_1,\ldots,y_{n-1}$;
$j\in \{1, \ldots,m_1\}$ such that $H\cap RN_{1j}\neq H\cap N_{1j}$.
Then from $(k,l)\geqslant (1,j)$ follows that $H\cap RN_{kl}\neq H\cap N_{kl}$.
\end{proposition}
\begin{proof}
It is clear that $H\cap N_{kl}\supseteq\gamma_l(H\cap N_{k1})$.
Let $\phi$ be the endomorphism of $F$ defined by $\phi(y_n)=1$, $\phi(y_j)=y_j$ for $j\neq n$.
Then $\phi(N_{kl})=\gamma_l(H\cap N_{k1})$, therefore $u=\phi(u)\in \gamma_l(H\cap N_{k1})$ for any $u\in H\cap N_{kl}$.
Thus $H\cap N_{kl}=\gamma_l(H\cap N_{k1})$.

We denote by $N$ the group $H\cap N_{11}$, by $\mathfrak{X}$ --- the fundamental ideal of ${\bf Z}[N]$.
Let $\{x_z | z\in I\}$ be a free set of generators of $N$, $\{\partial_z | z\in I\}$ --- the Fox derivatives of ${\bf Z}[N]$ $(I\subseteq {\bf N})$. Since free set of generators of $H$ contains more than one element and
$N$ --- a normal subgroup of $H$ it follows that free set of generators of $N$ contains more than one element. So we may assume that $|I| > 1$.

We now show that
\begin{eqnarray}\label{prp_1_f_1}
H\cap RN_{1i}> \gamma_i(N),\,i=j,\ldots,m_1+1.
\end{eqnarray}
If $j=m_1+1$ the truth of the statement is obvious.
Assume inductively that $H\cap RN_{1i}> \gamma_i(N)\,(i=j,\ldots, l;\,l\leqslant m_1)$.
We need to show that $H\cap RN_{1,l+1}> \gamma_{l+1}(N)$.
Let $v\in (H\cap RN_{1l})\setminus \gamma_l(N)$.
By Lemma~\ref{lm_3} we may assume that $\partial_1(v)\notin \mathfrak{X}^{l-1}$ (without loss of generality).
We denote by $w$ the element $[v,x_2]$.

From $v\in H\cap RN_{1l}$ follows that $w\in H\cap RN_{1,l+1}$.
By Lemma~\ref{lm_6}, $\partial_1(w)\notin \mathfrak{X}^l$, hence by Lemma~\ref{lm_3} $w\notin \gamma_{l+1}(N)$
and, by induction on $l$, $H\cap RN_{1l}> \gamma_l(N)$ for each term $N_{1l}\,(l\geqslant j)$ of series {\rm (\ref{tm4_0})}.
From {\rm (\ref{prp_1_f_1})} follows that $H\cap RN_{21}> H\cap N_{21}$.

Since $H\cap RN_{kl}\geqslant \gamma_l(H\cap RN_{k1})$ and $H\cap N_{kl}=\gamma_l(H\cap N_{k1})$
we reach the conclusion by noting that if $H\cap RN_{k1}> H\cap N_{k1}$ then
$\gamma_{l}(H\cap RN_{k1})> \gamma_{l}(H\cap N_{k1})$.
\end{proof}

\begin{proposition}\label{prp_2}
Suppose $F$ is a free group on free generators $y_1,\ldots,y_n$, $1\neq N_{11}$ --- a normal subgroup of $F$, $F/N_{11}$ --- right-ordered and relatively free group,
\begin{eqnarray}\label{prp_2_f_1}
N_{11} > \ldots > N_{1,m_1+1}=N_{21} > \ldots > N_{s,m_s+1},
\end{eqnarray}
where $N_{kl}$ --- the $l$-th term of the lower central series of $N_{k1}$.
Let $r$ be an element of $N_{1i}\backslash N_{1,i+1}\,(i\leqslant m_1)$, $R$ --- the normal subgroup of $F$ generated by $r$, $H$ --- the subgroup of $F$, generated by $y_1,\ldots,y_{n-1}$.
Then from $H\cap RN_{21}=H\cap N_{21}$ follows that
$H\cap RN_{kl}=H\cap N_{kl}$ for each term $N_{kl}\,(k> 1)$ of series {\rm (\ref{prp_2_f_1})}.
\end{proposition}
\begin{proof}
If $n=2$ then $H\cap N_{21}=1$. By the assumptions of the proposition $H\cap RN_{21}=H\cap N_{21}$ hence
$H\cap RN_{21}=1$ and $H\cap RN_{kl}=H\cap N_{kl}=1,\,k> 1$.
So we may assume that $n>2$.

We denote by $D_1,\ldots,D_n$ the Fox derivatives of ${\bf Z}[F]$ and denote by $\mathfrak{R}_{kl}$ the ideal ${\bf Z}[F]\cdot(\sqrt{RN_{kl}}-1)$.
Let us now suppose that $D_n(r)\equiv 0\mod {\mathfrak{R}_{21}}$.
Lemma~\ref{tm_3} tells us that there exist an elements $v_1,\ldots,v_d$ of $H\cap \sqrt{RN_{21}}$; $f_1,\ldots,f_d$ of $F$ such that $r\equiv v_1^{f_1}\cdots v_d^{f_d}\mod {[\sqrt{RN_{21}}\,,\sqrt{RN_{21}}\,]}$. Since $H\cap RN_{21}=H\cap N_{21}$ it follows that
$v_1,\ldots,v_d$ --- an elements of $N_{21}$.
Hence $r\in N_{1,i+1}$, a contradiction. So we may assume that $D_n(r)\not\equiv 0\mod {\mathfrak{R}_{21}}$.

It is clear that if $H\cap \sqrt{RN_{k,l+1}}=H\cap N_{k,l+1}$ then $H\cap RN_{k,l+1}=H\cap N_{k,l+1}$.
By assumption, $H\cap RN_{21}=H\cap N_{21}$, hence $H\cap \sqrt{RN_{21}}=H\cap RN_{21}=H\cap N_{21}$.
Now assume inductively that $H\cap \sqrt{RN_{ij}} = H\cap N_{ij}\,(i=2,\ldots, k;\,j=1,\ldots, l;\,l\leqslant m_k)$.
We need to show that $H\cap \sqrt{RN_{k,l+1}}=H\cap N_{k,l+1}$.

Let us denote the group $\sqrt{RN_{k1}}$ by $N$; the group $\sqrt{RN_{km}}$ by $N_m$ , $m\in {\bf N}$;
the fundamental ideal of ${\bf Z}[N]$ by $\mathfrak{X}$; the fundamental ideal of ${\bf Z}[H\cap N]$ by $\mathfrak{X^\prime}$.
Let $u\rightarrow \bar{u}$ be a Schreier coset representative function for $F\mod{N}$; $S$ --- a system of representatives of $F$ by $N$; $\{x_z | z\in {\bf N}\}$--- a free set of generators of $N$ such that $\{x_z | z\in {\bf N}\}\subseteq \{sy_j\overline{sy_j}^{\,-1} \mid s\in S,\, j=1,\ldots,n\}$; $\{\partial_z | z\in {\bf N}\}$ --- the Fox derivatives of ${\bf Z}[N]$;
$v\in H\cap N_{l+1}$.

From $N_{l+1}=\sqrt{RN_{k,l+1}}=\sqrt{R\gamma_{l+1} (N)}$ follows the existence of an element $u\in \gamma_{l+1} (N)$ and $j\in {\bf N}$ such that $v^j u^{-1}\in R$.
From (\ref{fm00}) follows the existence of an elements $k_p\in {\bf Z}[N]$, $f_p\in S$, $p=1,\ldots,d$ $(f_p\neq f_j\mbox{ for }p\neq j)$ such that $D_m(v^j u^{-1})\equiv D_m(r)\cdot  \sum_{p={1}}^{d} f_p k_p\mod{\mathfrak{R}_{kl}},\,m=1,\ldots,n$.

By (\ref{lm_5_f_1}), $D_m(u)=\sum_{z\in {\bf N}} D_m(x_z)\partial_z(u)$ thus
\begin{eqnarray}\label{prp_2_f_3}
D_m(v^j)\equiv D_m(r)\cdot  \sum_{p={1}}^{d} f_p k_p+ \sum_{z\in {\bf N}} D_m(x_z)\partial_z(u)\mod {\mathfrak{R}_{kl}},
\end{eqnarray}
$m=1,\ldots,n$. From $v\in H$ and (\ref{prp_2_f_3}) follows
\begin{eqnarray}\label{prp_2_f_4}
0\equiv D_n(r)\cdot  \sum_{p={1}}^{d} f_p k_p+ \sum_{z\in {\bf N}} D_n(x_z)\partial_z(u)\mod {\mathfrak{R}_{kl}}.
\end{eqnarray}
From $D_n(r)\not\equiv 0\mod {\mathfrak{R}_{k1}}$ follows the existence of an elements
$0\neq\gamma_p\in {\bf Z}$, $g_p\in S$, $p=1,\ldots,q;\,g_p\neq g_j$ for $p\neq j$ such that $D_n(r) \equiv  \sum_{p={1}}^q  \gamma_p g_p \mod {\mathfrak{R}_{k1}}$.

Let us now suppose that there exist $l_0,\,j$
such that $l_0< l$; $k_j\in \Delta_{l_0}\setminus \Delta_{l_0 +1}$; $k_p\in \Delta_{l_0}$, $p=1,\ldots,d$.
We may assume without loss of generality that $\{k_1,\ldots,k_e\}\subseteq \Delta_{l_0}\setminus\Delta_{l_0 +1}$
and $\{k_{e+1},\ldots,k_d\}\subseteq \Delta_{l_0 +1}$.
Since $F/N$ is a right-ordered group it follows that
there exists the element $g_a f_b$ of $\{g_pf_t\mid t\leqslant e\}$ such that $g_a f_b\neq g_c f_d\mod {N}$ for $(a,b)\neq (c,d)$ $(d\leqslant e)$.\\
Let $M\in S$ and $M\equiv g_a f_b\mod {N}$. We have
\begin{eqnarray}\label{prp_2_f_4-2}
D_n(r)\cdot  \sum_{p={1}}^{d} f_p k_p = M t_{ab} + h,
\end{eqnarray}
where $t_{ab}\in \Delta_{l_0}\setminus\Delta_{l_0 +1}$, $h$ --- a linear combination modulo $S\Delta_{l_0 +1}$ of an elements of the form $g t$, $t\in {\bf Z}[N]$, $g\in S$, $g\neq M$.

Since $\mathfrak{R}_{kl}=S\cdot{\bf Z}[N](N_l-1)\subseteq S\Delta_l$ and $\sum_{z\in {\bf N}} D_m(x_z)\partial_z(u)\in S\cdot \mathfrak{X}^l$ it follows that (\ref{prp_2_f_4-2}) contradicts (\ref{prp_2_f_4}). So we may assume that
\begin{eqnarray}\label{tm2_4_1_1}
k_p\in \Delta_l,\,p=1,\ldots,d.
\end{eqnarray}
By (\ref{prp_2_f_3}), (\ref{tm2_4_1_1}) we obtain
\begin{eqnarray}\label{prp_2_f_5}
D_m(v^j)= \sum_{p={1}}^{d_m} f_{pm} v_{pm},\,m=1,\ldots,n,
\end{eqnarray}
where $f_{pm}\in S$, $v_{pm}\in \Delta_l$.

If $x_z=Ky_m\overline{Ky_m}^{\,-1}$, $K\in S$ then Lemma~\ref{lm_7} tells us that
\begin{eqnarray}\label{prp_2_f_6}
D_m(v^j)= \overline{Ky_m}^{\,-1}\partial_z(v^j)+V,
\end{eqnarray}
where $V$ --- a linear combination of an elements of the form $g t$, $t\in {\bf Z}[N]$, $g\in S$, $g\not\equiv \overline{Ky_m}^{\,-1}\mod {N}$.

From (\ref{fm04}), (\ref{prp_2_f_5}), (\ref{prp_2_f_6}) it follows that $\partial_z(v^j)\in {\bf Z}[H_1]\cap\Delta_l=\Delta_l^\prime,\,z\in {\bf N}$.
We have
\begin{eqnarray}\label{prp_2_f_8}
H\cap N_t=H\cap N_{kt} = \gamma_t (H\cap N_{k1}), \,t=1,\ldots,l.
\end{eqnarray}
By $\gamma_t (H\cap N_{k1})-1\subseteq (\mathfrak{X}^\prime)^t$
and (\ref{prp_2_f_8}) we obtain $\Delta_l^\prime \subseteq (\mathfrak{X}^\prime)^l$, i.e.
$\partial_z(v^j)\in (\mathfrak{X}^\prime)^l,\,z\in {\bf N}$ whence
$v^j-1\in (\mathfrak{X}^\prime)^{l+1}$ and $v^j\in \gamma_{l+1} (H\cap N_{k1})$ \cite{Fx}.
Thus $v^j\in H\cap N_{k,l+1}$ hence $v\in H\cap N_{k,l+1}$ and, by induction on $l$,
$H\cap RN_{kl}=H\cap N_{kl}$ for each term $N_{kl}\,(k> 1)$ of series {\rm (\ref{prp_2_f_1})}.
\end{proof}

\begin{lemma}\label{lm_2}
Suppose $F$ is a free group on free generators $y_1,\ldots,y_n$; $H$ --- the subgroup of $F$, generated by $y_1,\ldots,y_{n-1}$; $N$ --- a normal subgroup of $F$, $F/N$ --- an orderable and relatively free group  with a free set $y_1N,\ldots,y_nN$ of generators;
$u\rightarrow \bar{u}$ --- a Schreier coset representative function for $F\mod{N}$, $S=S_\alpha\cup S_\beta$ --- a system of representatives of $F$ by $N$; $\delta_1,\ldots,\delta_l$, $\mu_1,\ldots,\mu_k$ --- an elements of $S$,
$\delta_iN <\delta_jN$, $\mu_iN <\mu_jN$ for $i <j$.
Then from $\{\overline{\mu_1^{-1}\mu_1},\ldots,\overline{\mu_1^{-1}\mu_k}\}\not\subseteq S_\alpha$ follows the existence
of an elements $\delta_{i_0},\,\mu_{j_0}$ such that $\overline{\delta_{i_0}\mu_{j_0}}\in S_\beta$ and $\overline{\delta_{i_0}\mu_{j_0}}\neq\overline{\delta_{i}\mu_{j}}$ for $(i_0,j_0)\neq (i,j)$.
\end{lemma}
\begin{proof}
Let $B$ be the normal subgroup of $F/N$ generated by $y_nN$, $A$ --- the subgroup of $F/N$ generated by $y_1N,\ldots,y_{n-1}N$. It is clear that $F/N=AB$, $A\cap B=1$.
We denote $\overline{\delta_i{\delta_1}^{-1}}N$ by $b_ia_i$ and $\overline{{\mu_1}^{-1}\mu_j}N$ by $\hat{a}_j\hat{b}_j$,
where $a_i,\hat{a}_j\in A$, $b_i,\hat{b}_j\in B$.

Let $\{\overline{\mu_1^{-1}\mu_1},\ldots,\overline{\mu_1^{-1}\mu_k}\}\not\subseteq S_\alpha$, $x=\max\,(b_1,\ldots,b_l)$,
$z=\max\,(\hat{b}_1,\ldots,\hat{b}_k)$, $b_{i_0}a_{i_0}=\max\,(b_ia_i\mid b_i=x)$, $\hat{a}_{j_0}\hat{b}_{j_0}=\max\,(\hat{a}_j\hat{b}_j\mid \hat{b}_j=z)$.
We may assume without loss of generality that $z>1$.

Since $\overline{\delta_1\mu_1}\neq\overline{\delta_{i}\mu_{j}}$ for $(1,1)\neq (i,j)$ we may assume that $\overline{\delta_1\mu_1}\in S_\alpha$. Then $t_{ij}=a_i(\overline{\delta_1\mu_1}N)\hat{a}_j\in A$. From $\overline{\delta_{i_0}\mu_{j_0}}N=b_{i_0}t_{i_0 j_0}\hat{b}_{j_0}$; ${\hat{b}_{j_0}}^{t_{i_0 j_0}}>1$;
$b_{i_0}\geqslant 1$ we obtain $\overline{\delta_{i_0}\mu_{j_0}}N=b_{i_0}{\hat{b}_{j_0}}^{t_{i_0 j_0}}t_{i_0 j_0}\notin A$.
If $\overline{\delta_{i_0}\mu_{j_0}}N = \overline{\delta_i\mu_j}N$ then
$t_{i_0 j_0}=t_{ij}$, $b_{i_0}=b_i$, $\hat{b}_{j_0}=\hat{b}_j$; $b_{i_0} a_{i_0}>b_i a_i$ for $i_0\neq i$; $\hat{a}_{j_0} \hat{b}_{j_0} > \hat{a}_j \hat{b}_j$ for $j_0\neq j$; $\overline{\delta_{i_0}\mu_{j_0}}N = b_{i_0} a_{i_0}(\overline{\delta_1\mu_1}N)\hat{a}_{j_0} \hat{b}_{j_0}>b_i a_i(\overline{\delta_1\mu_1}N)\hat{a}_j \hat{b}_j=\overline{\delta_i\mu_j}N$ for $(i_0,j_0)\neq (i,j)$. Hence $\overline{\delta_{i_0}\mu_{j_0}}N\neq \overline{\delta_i\mu_j}N$ for $(i_0,j_0)\neq (i,j)$.
\end{proof}

\begin{proposition}\label{prp_3}
Suppose $F$ is a free group on free generators $y_1,\ldots,y_n$,
$H$ --- the subgroup of $F$, generated by $y_1,\ldots,y_{n-1}$, $1\neq N$ --- a normal subgroup of $F$, $F/N$  --- an orderable and relatively free group,
\begin{eqnarray}\label{prp_3_f_1}
N=N_{11} > \ldots > N_{1,m_1+1}=N_{21} > \ldots > N_{s,m_s+1},
\end{eqnarray}
where $N_{kl}$ --- the $l$-th term of the lower central series of $N_{k1}$.
Let $r$ be an element of $N_{1i}\backslash N_{1,i+1}\,(i\leqslant m_1)$, $R$ --- the normal subgroup of $F$ generated by $r$.
Then from $r$ is not conjugate to any element of $HN_{1,i+1}$ follows that $H\cap RN_{1l}=H\cap N_{1l}\,(l=1,\ldots, m_1+1)$.
\end{proposition}
\begin{proof}
\noindent If $n=2$ and $N\neq F$ then $H\cap RN_{1l}=H\cap N_{1l}=1$.
If $n=2$ and $N= F$ then $H\cap RN_{11}=H\cap N_{11}=H$, $H\cap RN_{1l}=H\cap N_{1l}=1,\,l>1$.
So we may assume that $n>2$.

It is clear that if $H\cap \sqrt{RN_{1,l+1}}=H\cap N_{1,l+1}$ then $H\cap RN_{1,l+1}=H\cap N_{1,l+1}$.
We have $H\cap \sqrt{RN_{11}}=H\cap N_{11}$.
Now assume inductively that $H\cap \sqrt{RN_{1j}} = H\cap N_{1j}\,(j=1,\ldots, l;\,l\leqslant m_1)$.
We need to show that $H\cap \sqrt{RN_{1,l+1}}=H\cap N_{1,l+1}$.

Let us denote the group $\sqrt{RN_{1m}}$ by $N_m$ , $m\in {\bf N}$;
the fundamental ideal of ${\bf Z}[N]$ by $\mathfrak{X}$; the fundamental ideal of ${\bf Z}[H\cap N]$ by $\mathfrak{X^\prime}$.
Let $u\rightarrow \bar{u}$ be a Schreier coset representative function for $F\mod{N}$, $S=S_\alpha\cup S_\beta$ --- a system of representatives of $F$ by $N$; $D_1,\ldots,D_n$ --- the Fox derivatives of ${\bf Z}[F]$; $\mathfrak{R}_l={\bf Z}[F]\cdot (N_l-1)$;
$\{x_z \mid z\in P\}$ --- a free set of generators of $N$ such that $\{x_z \mid z\in P\}\subseteq \{sy_j\overline{sy_j}^{\,-1} \mid s\in S,\, j=1,\ldots,n\}$, where $P = {\bf N}$ for $N\neq F$ and $P = \{1,\ldots,n\}$ for $N=F$;
$\{\partial_z \mid z\in P\}$ ---  the Fox derivatives of ${\bf Z}[N]$; $v\in H\cap N_{l+1}$.

Free generator $x_z$ we name $\alpha$-generator if
$x_z\in  \{sy_j\overline{sy_j}^{\,-1} \mid s\in S_\alpha,\, j=1,\ldots,n-1\}$.
Thus if $x_z$ is not $\alpha$-generator then
$x_z\in  \{sy_n\overline{sy_n}^{\,-1} \mid s\in S\}$ or $x_z\in  \{sy_j\overline{sy_j}^{\,-1} \mid s\in S_\beta,\, j=1,\ldots,n-1\}$.

From $N_{l+1}=\sqrt{RN_{1,l+1}}=\sqrt{R\gamma_{l+1} (N)}$ follows the existence of an elements $u\in \gamma_{l+1} (N)$ and $j\in {\bf N}$ such that $v^j u^{-1}\in R$.
From (\ref{fm00}) follows the existence of an element $A\in {\bf Z}[F]$ such that $D_m(v^j u^{-1})
\equiv D_m(r)\cdot  A\mod{\mathfrak{R}_l},\,m=1,\ldots,n$.
By (\ref{lm_5_f_1}), $D_m(u)=\sum_{z\in {\bf N}} D_m(x_z)\partial_z(u)$ thus
\begin{eqnarray}\label{prp_3_f_2}
D_m(v^j)\equiv D_m(r)\cdot A + \sum_{z\in {\bf N}} D_m(x_z)\partial_z(u)\mod {\mathfrak{R}_l},
\end{eqnarray}
$m=1,\ldots,n$.  From $v\in H$ and (\ref{prp_3_f_2}) follows
\begin{eqnarray}\label{prp_3_f_3}
0\equiv D_n(r)\cdot A + \sum_{z\in {\bf N}} D_n(x_z)\partial_z(u)\mod {\mathfrak{R}_l}.
\end{eqnarray}
Since $v^j\in N_{1l}$, $r\in N_{1i}$, $u\in \gamma_{l+1} (N)$ it follows that $D_m(v^j)=\sum_{z\in {\bf N}} D_m(x_z)\partial_z(v^j)\in S\Delta_{l-1}$,
$D_m(r)\in S\Delta_{i-1}$,
$D_m(u)\in S\Delta_l$ $(m=1,\ldots,n)$.

Let us now suppose that there exists $m\in P$ such that $\partial_m(v^j)\not\in \Delta_l$.
If $x_m=Ky_t\overline{Ky_t}^{\,-1}$, $K\in S$, then Lemma~\ref{lm_7} tells us that
$D_t(v^j)= \overline{Ky_t}^{\,-1}\partial_m(v^j)+V$,
where $V$ --- a linear combination of an elements of the form $g t$, $g\in S$, $g\not\equiv \overline{Ky_t}^{\,-1}\mod {N}$, $t\in {\bf Z}[N]$. Hence $D_t(v^j)\in S\Delta_{l-1}\setminus S\Delta_l$.
Then from (\ref{prp_3_f_2}) follows that $D_t(r)\in S\Delta_{i-1}\setminus S\Delta_i$ and $A\in S\Delta_{l-i}\setminus S\Delta_{l-i+1}$. Since $D_t(v^j)\in S_\alpha {\bf Z}[N]$, by Lemma \ref{lm_2} there exists
$\mu\in S$ such that $\mu^{-1}A\in S_\alpha{\bf Z}[N]\mod{S\Delta_{l-i+1}}$.

We now show that if $x_z$ is not a $\alpha$-generator then
\begin{eqnarray}\label{prp_3_f_5}
\partial_z(r^\mu)\in\Delta_i.
\end{eqnarray}
Let $\partial_z(r^\mu)\in\Delta_{i-1}\setminus\Delta_i$, $x_z=Ky_t\overline{Ky_t}^{\,-1}$ where ($t=n$ and $K\in S$) or ($t<n$ and $K\in S_\beta$).
Then $D_t(r^\mu)\equiv \overline{Ky_t}^{\,-1}\partial_z(r^\mu) + V \mod{{\bf Z}[F]\cdot\Delta_i}$,
where $V$ --- a linear combination of an elements of the form $g t$, $g\in S$, $g\not\equiv \overline{Ky_t}^{\,-1}\mod {N}$,
$t\in \Delta_{i-1}\setminus \Delta_i$.
Thus we have $D_t(r^\mu)\cdot \mu^{-1} A \not\in{\bf Z}[F]\cdot\Delta_l$ which for $t=n$
contradicts (\ref{prp_3_f_3}) and for $t<n$, taking into account that $\overline{Ky_t}^{\,-1}\partial_z(r^\mu)\in S_\beta\Delta_{i-1}$,
contradicts (\ref{prp_3_f_2}).

Hence (\ref{prp_3_f_5}) holds. We now show that if $x_z=Ky_t\overline{Ky_t}^{\,-1}$, $t<n$, $K\in S_\alpha$, $\partial_z(r^\mu)\not\in\mathfrak{X}^i$ then
there exists $M\in {\bf N}$ such that
\begin{eqnarray}\label{prp_3_f_9}
M\partial_z(r^\mu)\in{\bf Z}[N\cap H]\mod{\Delta_i}.
\end{eqnarray}
We have
\begin{eqnarray}\label{prp_3_f_8}
D_t(r^\mu)\equiv \overline{Ky_t}^{\,-1}\partial_z(r^\mu) + L \mod{{\bf Z}[F]\cdot\Delta_i},
\end{eqnarray}
where $L$ --- a linear combination of an elements of the form
$g t$, $t\in \Delta_{i-1}\setminus \Delta_i$, $g\in S_\alpha$, $g\not\equiv \overline{Ky_t}^{\,-1}\mod {N}$.
From $D_t(v^j)\equiv D_t(r^\mu)\cdot  (\mu^{-1}A)\mod{{\bf Z}[F]\cdot\Delta_l}$
and (\ref{prp_3_f_8}) by Lemma \ref{lm_1} we obtain (\ref{prp_3_f_9}).
By Lemma \ref{gr1_lm2} from (\ref{prp_3_f_5}), (\ref{prp_3_f_9}) follows that $r^{M\mu}\in HN_{1,i+1}$.
Since $N/N_{1,i+1}$ is a nilpotent torsion free group
it follows that if $x,\, y$ --- an elements of $N/N_{1,i+1}$ and $x^n=y^n\,(n\neq 0)$ then $x=y$. Thus $r^\mu\in HN_{1,i+1}$, a contradiction.

So we may assume that $\{\partial_z(v^j) \mid z\in P\} \subseteq \Delta_l$.
From $\partial_z(v^j)\in {\bf Z}[H_1]\cap\Delta_l$ by (\ref{fm04}) we obtain $\partial_z(v^j)\in \Delta_l^\prime,\,z\in P$.
If $t\leqslant l$ then $H\cap N_t=H\cap N_{1t}=\gamma_t (H\cap N)$. Hence
$\Delta_l^\prime\subseteq (\mathfrak{X}^\prime)^l$ and
$\partial_z(v^j)\in (\mathfrak{X}^\prime)^l,\,z\in P$;
$v^j-1\in (\mathfrak{X}^\prime)^{l+1}$; $v^j\in \gamma_{l+1} (H\cap N)$ \cite{Fx}.
Thus $v^j\in H\cap N_{1,l+1}$ whence $v\in H\cap N_{1,l+1}$ and
$H\cap RN_{1l}=H\cap N_{1l}$ for each term $N_{1l}$ of series {\rm (\ref{prp_3_f_1})}.
\end{proof}
From Propositions \ref{prp_1}, \ref{prp_2}, \ref{prp_3}  we obtain immediately Theorem~\ref{tm_1}.

\section{A generalized theorem on freedom for relatively free groups}

Suppose $G$ is a group. The elementary transformations
of a matrix over ${\bf Z}[G]$ are defined as follows:
\begin{align}
&\text{Interchange of the columns $i$ and $j$;}\label{df_1}\\
&\text{Interchange of the rows $i$ and $j$;}\label{df_2}\\
&\text{Right-multiply the }i\text{-th row by a non-zero element of }{\bf Z}[G];\label{df_3}\\
&\text{Add to the }j\text{-th row the }i\text{-th row, }\label{df_4}\\
&\text{right-multiplied by a non-zero element of }{\bf Z}[G]\text{, where }i<j.\notag
\end{align}

Let $M=\|m_{kn}\|$ be a $r\times s$ matrix over ${\bf Z}[G]$, $t$ --- the rank of $M$. We denote by $\Phi(M)$ the matrix obtained from $M$ by a chain $\Phi$ of elementary transformations of $M$; $M$ is said to be the lower triangular matrix if $m_{kk}\neq 0\,(k=1,\ldots, t)$, $m_{kn}= 0\,(k>n\,or\, k>t)$.\\
Let $v$ be a $1\times s$ matrix over ${\bf Z}[G]$. We define $\Phi(v)$
by putting $\psi(v)=v$ for any elementary transformation of rows $\psi$ of $\Phi$.

\begin{lemma}\label{lm_8}
Let $F$ be a free group, $N=N_1 \geqslant \ldots \geqslant N_m\geqslant \ldots$  --- a descending series of normal subgroups
of $F$ with abelian torsion free factors, $[N_p\,,N_q\,]\leqslant N_{p+q}$,
$F/N$ --- a soluble right-ordered group, $\phi$ --- natural homomorphism ${\bf Z}[F]\to {\bf Z}[F/N_m]$.\\
Let $\|a_{kn}\|$ be a $r\times s$ matrix over ${\bf Z}[F/N_m]$, $\phi^\prime$ --- natural homomorphism ${\bf Z}[F/N_m]\to {\bf Z}[F/N]$,  $\phi^\prime(a_{kk})\neq 0$, if $n<k$, then $\phi^\prime(a_{kn})= 0$ $(k=1,\ldots, r)$, $\psi$ --- a valuation on ${\bf Z}[F/N_m]$ defined by $\psi(\phi(u))=j$ if $u\in S\Delta_j\setminus S\Delta_{j+1}\mod{{\bf Z}[F](N_m-1)}$, $\psi(0)=\infty$.
Then matrix $\|a_{kn}\|$ by a finite number of operations {\rm (\ref{df_3})}, {\rm (\ref{df_4})}
can be converted into a matrix $\|b_{kn}\|$ such that $\psi(b_{kk})\leqslant \psi(b_{kn})$; $b_{kk}\neq 0$; if $n<k$, then $b_{kn}= 0$ (k=1,\ldots, r; n=1,\ldots, s).
\end{lemma}
\begin{proof}
Since $F/N_m$ is a soluble torsion free group, ${\bf Z}[F/N_m]$ is known to satisfy right Ore's condition \cite{Lv}.
The groups $F/N$, $F/N_m$ are right-ordered, hence ${\bf Z}[F/N]$, ${\bf Z}[F/N_m]$ have no zero-divisors.

By hypothesis, $\phi^\prime(a_{kk})\neq 0$, if $n<k$, then $\phi^\prime(a_{kn})= 0$.
Hence $\psi(a_{kk})=0$ and if $n<k$ then $\psi(a_{kn})>0$, $k=1,\ldots,r$.

We put $(b_{11},\ldots,b_{1s})=(a_{11},\ldots,a_{1s})$. It is clear that $\psi(b_{11})\leqslant \psi(b_{1j})$, $j=1,\ldots,s$.
Now assume inductively that by a finite number of operations {\rm (\ref{df_3})}, {\rm (\ref{df_4})}
the rows $(a_{k1},\ldots,a_{ks})$ can be converted into $(b_{k1},\ldots,b_{ks})$
such that $b_{kk}\neq 0$; if $n<k$ then $b_{kn}= 0$; $\psi(b_{kk})\leqslant \psi(b_{kn})$, $k=1,\ldots, t-1; n=1,\ldots,s$.

If $a_{t1}=\ldots=a_{t,t-1}=0$ then we put $(b_{t1},\ldots,b_{ts})=(a_{t1},\ldots,a_{ts})$. It is clear that
$b_{tt}\neq 0$; if $n<t$ then $b_{tn}= 0$; $\psi(b_{tt})\leqslant \psi(b_{tj})$, $j=1,\ldots,s$.

Let us now suppose that $a_{t1}=\ldots=a_{t,l-1}=0$, $a_{t,l}\neq 0$, $l\leqslant t-1$. There exist a non-zero elements $\beta_1,\,\beta_2$ of ${\bf Z}[F/N_m]$ such that $b_{ll}\beta_1=-a_{tl}\beta_2$.\\
We put $c_{tn}=b_{ln}\beta_1+a_{tn}\beta_2$, $n=1,\ldots,s$. Then $c_{t1}=\ldots=c_{tl}= 0$.

From $\psi(b_{ll}\beta_1)\leqslant \psi(b_{lj}\beta_1)$ $(j=1,\ldots,s)$, $\psi(b_{ll}\beta_1)=\psi(a_{tl}\beta_2)$
we obtain\\
$\psi(a_{tt}\beta_2)< \psi(a_{tl}\beta_2)\leqslant \psi(b_{lj}\beta_1)$ $(j=1,\ldots,s)$.\\
Since $\psi(a_{tt}\beta_2)< \psi(a_{tj}\beta_2)$ $(j<t)$ and $\psi(a_{tt}\beta_2)\leqslant \psi(a_{tj}\beta_2)$ $(j=1,\ldots,s)$,
we have $\psi(c_{tt})=\psi(a_{tt}\beta_2)< \psi(c_{tj})$ $(j<t)$ and $\psi(c_{tt})\leqslant \psi(c_{tj})$, $j=1,\ldots,s$.

Thus the row $(a_{t1},\ldots,a_{ts})$ can be converted by a finite number of operations {\rm (\ref{df_3})}, {\rm (\ref{df_4})} into $(b_{t1},\ldots,b_{ts})$ such that $b_{tt}\neq 0$; if $n<t$ then $b_{tn}= 0$; $\psi(b_{tt})\leqslant \psi(b_{tj})$ $(j=1,\ldots,s)$ and the proof is complete.
\end{proof}

Let $\|a_{kn}\|$ be a $r\times s$ matrix over ${\bf Z}[F/N_m]$, $t$ --- the rank of $\|a_{kn}\|$.
It is not hard to verify that $\|a_{kn}\|$ can be converted by a finite number of operations
{\rm (\ref{df_1})}-{\rm (\ref{df_4})}
into a lower triangular form $\|b_{kn}\|$ such that $\psi(b_{kk})\leqslant \psi(b_{kn})\,(k=1,\ldots, t; n=1,\ldots, s)$.

\begin{lemma}\label{lm_9}
Suppose $G$ is a soluble torsion free group; $M$ --- $r\times s$ matrix over ${\bf Z}[G]$; $\alpha_i$ --- $i$-th row of $M$;
$\alpha$ --- a right-linear combination of rows $\alpha_1,\ldots,\alpha_r$.
Then for each elementary transformation $\psi$ of $M$ there exists non-zero $d_\psi\in {\bf Z}[G]$ such that $\psi(\alpha) d_\psi$ is a right-linear combination of rows of $\psi(M)$.
\end{lemma}
\begin{proof}
If $\psi$ is one of the operations  (\ref{df_1}),
(\ref{df_2}), (\ref{df_4}) the result is obvious.
Let $\psi(M)$ be a matrix obtained by operation (\ref{df_3}); $a$ ---
non-zero element of ${\bf Z}[G]$; $\alpha_1,\ldots,\alpha_i a,\ldots,\alpha_r$ --- the rows of
$\psi(M)$.
By assumption, there exist an elements $b_1,\ldots,b_r$ of ${\bf Z}[G]$ such that
$\alpha_1b_1+\ldots+\alpha_ib_i+\ldots+\alpha_rb_r=\alpha$. It is no restriction to
assume that $b_i\neq 0$, otherwise the result is obvious.
Since $G$ is a soluble torsion free group, ${\bf Z}[G]$ is known to satisfy right Ore's condition \cite{Lv}.
Thus there exist non-zero $c$, $d_\psi$ such that $ac=b_id_\psi$ and we have $\alpha_1b_1d_\psi+\ldots+\alpha_i ac+\ldots+\alpha_rb_rd_\psi=\alpha\, d_\psi$.
\end{proof}

The proof of Theorem \ref{tm_2}.
We may clearly assume that $n-m>1$.
We denote by $\phi_k$ the natural homomorphism
${\bf Z}[F]\to {\bf Z}[F/\sqrt{RN_{k,m_k+1}}\,]$; by $\phi_k^\prime$ --- natural homomorphism ${\bf Z}[F/\sqrt{RN_{k,m_k+1}}\,]\to {\bf Z}[F/\sqrt{RN_{k1}}\,]$; by $\phi_0$ --- natural homomorphism ${\bf Z}[F]\to {\bf Z}[F/N_{11}]$; by $\phi^\prime_0$ --- natural homomorphism ${\bf Z}[F/N_{11}]\to {\bf Z}[F/N_{11}]$. Let $A=\|a_{rt}\|$  be a matrix over ${\bf Z}[F]$. Then we denote by $A^{\phi_k}$ --- the matrix $\|\phi_k(a_{rt})\|$ and denote by $(A^{\phi_k})^{\phi_k^\prime}$ the matrix $\|\phi_k^\prime(\phi_k(a_{rt}))\|$.

We denote by $D_1,\ldots,D_n$ the Fox derivations of the group ring ${\bf Z}[F]$;
by $m_{ij}$ the elements $D_j(r_i)\,(i=1,\ldots, m;\,j=1,\ldots, n)$; by $M$ the matrix $\|m_{ij}\|$; by $t_k$ --- the rank of $M^{\phi_k}$.

Let $R\subseteq N_{k,l}$ $(l<m_k+1)$. If $m_{ij}\equiv 0 \mod{{\bf Z}[F](\sqrt{RN_{s,m_s+1}}-1)}$ for any $m_{ij}$ then $m_{ij}\equiv 0 \mod{{\bf Z}[F](N_{k,l}-1)}$ for any $m_{ij}$ whence $r_i\in \gamma_2(N_{k,l})\subseteq N_{k,l+1}\,(i=1,\ldots, m)$ \cite{Sch}. Thus if $m_{ij}\equiv 0 \mod{{\bf Z}[F](\sqrt{RN_{s,m_s+1}}-1)}$ for any $m_{ij}$
then $R\subseteq N_{s,m_s+1}$ and $H\cap RN_{kl} = H\cap N_{kl}$ for each term $N_{kl}$ of series {\rm (\ref{tm_2_f_1})}, where $H=F$.

Further we may assume that there exists $K\in \{0,\ldots,s\}$ such that $t_K>0$ and if $i<K$ then $t_i=0$.
Let $\psi_k$  be the valuation on ${\bf Z}[F/\sqrt{RN_{km}}]$ defined by, $\Phi_K$ --- a chain of elementary transformations of $M$ such that $M_K=\|m^{(K)}_{ij}\|=(\Phi_K(M))^{\phi_K}$ be a lower triangular $m\times n$ matrix and $\psi_K(m^{(K)}_{ii})\leqslant \psi_K(m^{(K)}_{ij})\,(i=1,\ldots, t_K;\,j=1,\ldots, n)$.

Let $K<s$. Now assume inductively that for some $k$ $(K\leqslant k< s)$ we have $\Phi_k$ --- a chain of elementary transformations of $M$; $M_k=\|m^{(k)}_{ij}\|=(\Phi_k(M))^{\phi_k}$ --- a lower triangular $m\times n$ matrix;
$\psi_k(m^{(k)}_{ii})\leqslant \psi_k(m^{(k)}_{ij})\,(i=1,\ldots, t_k;\,j=1,\ldots, n)$.

We denote $\Phi_k(M)$ by $M_{k+1,1}$. We have $(M_{k+1,1}^{\phi_{k+1}})^{\phi^\prime_{k+1,1}}=M_k$. Then Lemma~\ref{lm_8} tells us that there exists a chain $\Phi_{k+1,1}$ of elementary transformations (\ref{df_3}) (where $i\leqslant t_k$), (\ref{df_4}) (where $i<j\leqslant t_k$) of $M_{k+1,1}$ such that $(\Phi_{k+1,1}(M_{k+1,1}))^{\phi_{k+1}}$ be $m\times n$ matrix $\|b_{ij}\|$ with $b_{ii}\neq 0$; $b_{ij}=0\,(j< i)$; $\psi_{k+1}(b_{ii})\leqslant \psi_{k+1}(b_{ij})\,(i=1,\ldots, t_k;\,j=1,\ldots, n)$.

We denote $\Phi_{k+1,1}(M_{k+1,1})$ by $M_{k+1,2}$
and denote by $\Phi_{k+1,2}$ a chain of operations (\ref{df_3}) (where $i> t_k$), (\ref{df_4}) (where $i\leqslant t_k\text{ and }j>t_k$) such that $(\Phi_{k+1,2}(M_{k+1,2}))^{\phi_{k+1}}$ be a matrix $\|c_{ij}\|$ with $c_{ij}=b_{ij}\,(i=1,\ldots, t_k;\,j=1,\ldots, n)$ and $c_{ij}=0\,(j\leqslant t_k<i)$.

We denote $\Phi_{k+1,2}(M_{k+1,2})$ by $M_{k+1,3}$
and denote by $\Phi_{k+1,3}$ a chain of operations (\ref{df_1})-(\ref{df_4}) (where $i,\,j> t_k$) of $M_{k+1,3}$ such that $(\Phi_{k+1,3}(M_{k+1,3}))^{\phi_{k+1}}$ be a lower triangular $m\times n$ matrix $\|m^{(k+1)}_{ij}\|$ and $\psi_{k+1}(m^{(k+1)}_{ii})\leqslant \psi_{k+1}(m^{(k+1)}_{ij})\,(i=1,\ldots, t_{k+1};$ $j=1,\ldots, n)$.
We denote the matrix $\|m^{(k+1)}_{ij}\|$ by $M_{k+1}$, the sequence $\Phi_k$, $\Phi_{k+1,1}$,
$\Phi_{k+1,2}$, $\Phi_{k+1,3}$ by $\Phi_{k+1}$.
Thus by induction on $k$ we have $M_k$, $\Phi_k$ for any $k\in \{K,\ldots,s\}$.
If $K=s$ the truth of the statement is obvious.

Let $I_s=\{i_1,\ldots,i_{t_s}\}$ be the subset of $\{1,\ldots,n\}$ such that
if $m_{i_j}$ --- $i_j$-th column of $M$ then $\Phi_s(m_{i_j})$ --- $j$-th column of $\Phi_s(M)$;
$\{j_1,\ldots,j_p\}= \{1,\ldots,n\}\setminus I_s$;
$H$ --- the free group with the free set  $\{y_{j_1},\ldots,y_{j_p}\}$ of generators.
Since $t_s\leqslant m$ it follows that $p\geqslant n-m$.

We denote by $N$ the group $\sqrt{RN_{k1}}$; by $\mathfrak{X}$ --- the fundamental ideal of ${\bf Z}[N]$; by $\mathfrak{X^\prime}$ --- the fundamental ideal of ${\bf Z}[H\cap N]$; by $N_m$ the group $\sqrt{RN_{km}}$.
Let $S$ be a system of representatives of $F$ by $N$. We define a valuation $\psi_k$ on ${\bf Z}[F/N_{m_k+1}]$
by putting $\psi_k(0)=\infty$ and $\psi_k(\phi_k(u))=j$ for any $u\in S\Delta_j$, $u\not\in S\Delta_{j+1}+{\bf Z}[F](N_{m_k+1}-1)$.
Let $\{x_{kz} | z\in P\}$ be a free set of generators of $N$, $\{\partial_{kz} | z\in P\}$ ---  the Fox derivatives of ${\bf Z}[N]$ $(P\subseteq {\bf N})$, $\mathfrak{R}_l={\bf Z}[F]\cdot (N_l-1)$,
$v\in H\cap (R\gamma_{m_k+1} (N))$.

From $v\in R\gamma_{m_k+1} (N)$ follows the existence of an element $u\in \gamma_{m_k+1} (N)$ such
that $vu^{-1}\in R$.
From (\ref{fm00}) follows the existence of an elements $\beta_1,\ldots,\beta_m\in {\bf Z}[F]$ such that
$D_j(vu^{-1})\equiv \sum_{i={1}}^{m} D_j(r_i)\beta_i\mod {\mathfrak{R}_{m_k+1}}$.
Hence
\begin{eqnarray}\label{tm_2_f_2}
D_j(v)\equiv \sum_{i={1}}^{m} D_j(r_i)\beta_i+D_j(u)\mod {\mathfrak{R}_{m_k+1}},\,j= 1,\ldots,n.
\end{eqnarray}
Since $\partial_{kz}(u)\in \mathfrak{X}^{m_k}$, $z\in P$ it follows that
\begin{eqnarray}\label{tm_2_f_2_2}
D_j(u)\in {\bf Z}[F]\cdot\mathfrak{X}^{m_k},\,j= 1,\ldots,n.
\end{eqnarray}
We now show that
\begin{eqnarray}\label{tm_2_f_3}
D_j(v)\in {\bf Z}[F]\cdot\Delta_{m_k},\,j= 1,\ldots,n.
\end{eqnarray}
If $t_k=0$ then $D_j(r_i)\equiv 0\mod {\mathfrak{R}_{m_k+1}}$, hence from
(\ref{tm_2_f_2}), (\ref{tm_2_f_2_2}) follows (\ref{tm_2_f_3}).
Let $t_k>0$, $V=(D_1(v)-D_1(u),\ldots,D_n(v)-D_n(u))$.
From Lemma~\ref{lm_9} and (\ref{tm_2_f_2}) we obtain the existence of an element $d$ of ${\bf Z}[F]$ such that $\phi_k(d)\neq 0$ and $(\Phi_k(V d))^{\phi_k}$ --- a right-linear combination of non-zero rows of $M_k$.

Since $v\in H$ it follows that $D_j(v)=0$ $(j=i_1,\ldots,i_{t_s})$ and therefore an element of $\Phi_k(V d)$ in the $i$-th coordinate belongs to
${\bf Z}[F]\Delta_{m_k}d\, (i=1,\ldots, t_k)$.
Since $\psi_k(m^{(k)}_{ii})\leqslant \psi_k(m^{(k)}_{ij})\,(i=1,\ldots, t_k;\,j=1,\ldots, n)$ it is not hard to verify that all elements of $(\Phi_k(V d))^{\phi_k}$ belong to $\phi_k({\bf Z}[F]\Delta_{m_k}d)$ and (\ref{tm_2_f_3}) holds.

Thus $D_j(v)-D_j(u)\in {\bf Z}[F]\cdot\Delta_{m_k}\mod{\mathfrak{R}_{m_k+1}}\,(j=1,\ldots,n)$, whence $D_j(v)\in {\bf Z}[F]\cdot\Delta_{m_k},\,j= 1,\ldots,n,$ as we wished to show.

It is clear that $H\cap \sqrt{RN_{11}} = H\cap N_{11}$.
Now assume inductively that $H\cap \sqrt{RN_{ij}} = H\cap N_{ij}\,(i=1,\ldots, k;\,j=1,\ldots, l;\,l\leqslant m_k)$.
We need to show that $H\cap \sqrt{RN_{k,l+1}} = H\cap N_{k,l+1}$.
Let $v\in H\cap (R\gamma_{l+1} (\sqrt{RN_{k1}}))$.

Consider case $l=1$.
We have $H\cap N=H\cap N_{k1}$ and $v\in H\cap (R\gamma_2 (N))$.
From $v\in R\gamma_2 (N)$ follows the existence of an element $u$ of $\gamma_2 (N)$ such
that $vu^{-1}\in R$.
Then from (\ref{fm00}) follows the existence of an elements $B_1,\ldots,B_m\in {\bf Z}[F]$ such
that
\begin{eqnarray}\label{tm2_6_gr_ends_2}
D_j(vu^{-1})\equiv \sum_{i={1}}^{m} D_j(r_i)B_i\mod {\mathfrak{R}_1},\,j= 1,\ldots,n.
\end{eqnarray}
Let $V=(D_1(v)-D_1(u),\ldots,D_n(v)-D_n(u))$.
Since $v\in H$ it follows that $D_j(v)=0$ $(j=i_1,\ldots,i_{t_s})$. From $u\in \gamma_2 (N)$ we have $D_j(u)=0\mod {\mathfrak{R}_1},\,j= 1,\ldots,n$.
From Lemma~\ref{lm_9} and (\ref{tm2_6_gr_ends_2}) we obtain the existence of an element $d$ of ${\bf Z}[F]$ such that
$\phi_{k-1}(d)\neq 0$ and $(\Phi_{k-1}(V d))^{\phi_{k-1}}$ --- a right-linear combination of rows of a lower triangular matrix $M_{k-1}$. Thus $(\Phi_{k-1}(V d))^{\phi_{k-1}}$ --- trivial row whence
$D_j(v)\equiv 0 \mod{{\bf Z}[F]\cdot(N_{k1}-1)}$ $(j=1,\ldots,n)$ and $v\in N_{k2}$ \cite{Sch}.

Consider case $l>1$.
Since free set of generators of $H$ contains more than one element and
$H\cap N$ --- a normal subgroup of $H$ it follows that free set of generators of $H\cap N$ contains more than one element.
Let us denote a free set $\{x_1,\,x_2,\ldots\}$ of generators of $H\cap N$
by $X$. Since $H\cap N_2-1=H\cap N_{k2}-1=\gamma_2 (H\cap N_{k1})-1\subseteq \mathfrak{X^\prime}^2$ it follows that
${\bf Z}[H]\cap \Delta_2=\Delta_2^\prime=\mathfrak{X^\prime}^2$.
If $x\in X$ then $x-1\notin \mathfrak{X^\prime}^2$,
whence $x-1\notin \Delta_2$.
From $v\in \gamma_l (H\cap N_{k1})$ follows that $\partial_z(v)\in \mathfrak{X}^{l-1},\,x_z\in X$. We now show that
\begin{eqnarray}\label{tm_2_f_4}
\partial_z(v)\in \Delta_l,\, x_z\in X.
\end{eqnarray}

Suppose that there exists $i$ such that $\partial_i(v)\in \Delta_{l-1}\setminus \Delta_l$.
We define elements $v_1,\,v_2,\ldots$ of $H\cap N$ inductively as follows: $v_1=[v,x_t]$, $v_{j+1}=[v_j,x_t]$ and
denote by $\bar{v}$ the element $v_{m_k -l}$ $(t\neq i)$. It is clear that $v_j\in \gamma_{l+j} (H\cap N_{k1})$ whence
$v_j - 1\in \mathfrak{X}^{l+j}$.
Since $\partial_i(v_{j+1})=-\partial_i(v_j)v_j^{-1}x_t^{-1}v_jx_t+\partial_i(v_j)x_t=-\partial_i(v_j)v_{j+1}+\partial_i(v_j)x_t$,
it follows that $\partial_i(v_{j+1})\equiv\partial_i(v_j)(x_t-1)\mod{\mathfrak{X}^{l+j+1}}$. Thus we have
\begin{eqnarray}\label{tm_2_f_5}
\partial_i(\bar{v}))\equiv\partial_i(v)(x_t-1)^{m_k -l}\mod{\mathfrak{X}^{m_k}}.
\end{eqnarray}
From $x_t-1\not\in \Delta_2$ and (\ref{tm_2_f_5}) follows that $\partial_i(\bar{v})\not\in \Delta_{m_k}$.
By Lemma \ref{lm_7}, there exist $D_j$, $f\in F$ such that
$D_j(\bar{v}) = f\partial_i(\bar{v})+V$,
where $V$ --- a linear combination of an elements of the form $g t$, $t\in {\bf Z}[N]$, $g\in S$, $g\notin f N$.
Thus we have $D_j(\bar{v})\not\in {\bf Z}[F]\cdot \Delta_{m_k}$
which contradicts (\ref{tm_2_f_3}).
This proves (\ref{tm_2_f_4}).

So we have $\partial_z(v)\in {\bf Z}[H_1]\cap \Delta_l$ and by (\ref{fm04}) $\partial_z(v)\in \Delta_l^\prime$.
If $t\leqslant l$ then $H\cap N_t = \gamma_t (H\cap N_{k1})$, hence
$\Delta_l^\prime\subseteq (\mathfrak{X}^\prime)^l$. We obtain that
$\partial_z(v)\in (\mathfrak{X}^\prime)^l$ whence
$v-1\in \mathfrak{X^\prime}^{l+1}$ and $v\in \gamma_{l+1} (H\cap N_{k1})$ \cite{Fx}. Thus from $v\in H\cap (R\gamma_{l+1} (\sqrt{RN_{k1}}))$ follows that $v\in H\cap N_{k,l+1}$.

From $v\in H\cap \sqrt{RN_{k,l+1}}$ follows that $v^c\in H\cap (R\gamma_{l+1} (\sqrt{RN_{k1}}))\text{ for some }c\neq 0$, whence
$v^c\in H\cap N_{k,l+1}$. Hence $v\in H\cap N_{k,l+1}$ and
$H\cap RN_{kl} = H\cap N_{kl}$ for each term $N_{kl}$ of series {\rm (\ref{tm_2_f_1})}.

\end{document}